\input amstex

\magnification=\magstep1  
\documentstyle{amsppt}
\hsize = 6.5 truein
\vsize = 9 truein
\NoRunningHeads

\topmatter
\title
Pairs of periodic orbits with fixed homology difference
\endtitle
\author
Morten S. Risager and Richard Sharp
\endauthor
\affil
University of Aarhus and University of Manchester
\endaffil
\address
Morten S. Risager, Department of Mathematical Sciences,
University of Aarhus,
Ny Munkegade Building 530,
8000 Aarhus,
Denmark.
\endaddress
\email
risager\@imf.au.dk
\endemail
\address
Richard Sharp, School of Mathematics, University
of Manchester, Oxford Road, Manchester M13 9PL, UK.
\endaddress
\email
sharp\@maths.man.ac.uk
\endemail
\abstract
We obtain an asymptotic formula for the number of pairs of closed
orbits of a weak-mixing transitive Anosov flow whose homology
classes have a fixed difference.
\endabstract
\endtopmatter
\document

\heading
0. Introduction
\endheading

Consider $M$ a compact smooth Riemannian manifold and 
$\phi_t : M \to M$ a transitive Anosov flow on $M$.
Such a manifold has a
countable infinity of (prime) periodic orbits $\gamma$. We denote the lenght af such an orbit by $l(\gamma)$. 
Writing $\pi(T) :=\#\{\gamma \hbox{ : } l(\gamma) \leq T\}$ the
following expansion holds: 
$\pi(T) \sim e^{hT}/hT$, as $T \to +\infty,$ 
where
$h>0$ is the topological entropy of $\phi$
\cite{11},\cite{12},\cite{14},\cite{15}.

To refine the problem one might try to understand the distribution
of periodic orbits with respect to the homology of $M$.
To keep our statements simple,
we shall suppose that $H_1(M,\Bbb Z)$ is infinite and ignore
any torsion.
Suppose that $M$ has first Betti number $k \geq 1$;
then we may fix an identification of  $H_1(M,\Bbb Z)/\text{torsion}$ 
with $\Bbb Z^k$.
For $\alpha \in \Bbb Z^k$, write
$\pi(T,\alpha) := \#\{\gamma \hbox{ : } l(\gamma) \leq T, \
[\gamma] = \alpha\},$
where $[\gamma]$ denotes the homology class of $\gamma$ (modulo torsion).
A variety of behaviours for this counting function are possible.
For example, for a geodesic flow (in variable negative curvature)
there exists $C>0$ (independent of $\alpha$) such that
$$\pi(T,\alpha) \sim C \frac{e^{hT}}{T^{1+k/2}}, 
\hbox{ as } T \to +\infty \eqno(1)$$ 
\cite{5},\cite{6},\cite{10},\cite{17},\cite{18}
but, for other Anosov flows, $\pi(T,\alpha)$ may grow at a slower rate
or even be bounded (or identically zero), depending on the
circumstances
\cite{7},\cite{20}.
Nevertheless, if $\alpha$ is allowed to grow 
with $T$ at an appropriate linear rate
then an asymptotic of the form (1) always holds \cite{2},\cite{9}.

In this note, we are going to study the \it relative distribution of
pairs \rm of closed orbits in $H_1(M,\Bbb Z)$.
For
$\beta \in \Bbb Z^k$, define
$$\pi_2^\beta(T) :=\#\{(\gamma,\gamma') \hbox{ : } l(\gamma),l(\gamma')
\leq T, \ [\gamma]-[\gamma']=\beta\}.
\eqno(2)$$
Since the asymptotic behaviour of $\pi(T,\alpha)$ is different for
different types of weak-mixing transitive Anosov flows one might
suspect that $\pi_2^\beta(T)$ has varying asymptotic behaviour as
well. We show that this is \it not \rm the case. Surprisingly the asymptotic
behaviour is universal. Our main result is the following:

\proclaim{Theorem 1}
Let $\phi_t : M \to M$ be a weak-mixing transitive Anosov flow
on a compact smooth Riemannian manifold $M$ with first Betti
number $k \geq 1$.
Then there exists $\Cal C(\phi)>0$ such that, for each $\beta \in 
H_1(M,\Bbb Z)/\text{\rom{torsion}} \cong \Bbb Z^k$,
$$\pi_2^\beta(T) \sim \Cal C(\phi) \frac{e^{2hT}}{T^{2+k/2}}, \hbox{
\text{as} }
T \to +\infty.$$
\endproclaim

\demo{Remark}
The constant $\Cal C(\phi)$ can be described in terms of the Hessian of 
an associated
entropy function at a special point. To be precise
$$\Cal C(\phi)=\frac{1}{2^k\pi^{k/2}\sigma^kh^2}$$ 
where $\sigma^{2k}$ is the determinant of minus the Hessian of this
entropy function evaluated at the winding cycle associated to the measure of
maximal entropy for $\phi$. See section 1 for details.
\enddemo

For a compact hyperbolic surface $V$ of genus $g \geq 2$, 
the geodesic flow on the sphere
bundle $SV$ is a weak-mixing transitive Anosov flow with topological
entropy equal to $1$.
Furthermore, the natural projection $p : SV \to V$ induces
an isomorphism between $H_1(SV,\Bbb Z)/\text{torsion}$
and $H_1(V,\Bbb Z) \cong \Bbb Z^{2g}$.   
There is a one-to-one correspondence between periodic
orbits for the flow and closed geodesics on the surface,
which preserves lengths and respects this isomorphism.
Thus $\pi_2^\beta(T)$
counts also the number of pairs of closed geodesics on $V$, with lengths
at most $T$, and with the two homology classes differing by
$\beta$.
We may recover the following result, previously
obtained by the first author \cite{19}.

\proclaim{Theorem 2}
Let $V$ be a compact hyperbolic surface of genus $g$.
Then, for each $\beta \in H_1(V,\Bbb Z) \cong \Bbb Z^{2g}$,
$$\pi_2^\beta(T) \sim \frac{(g -1)^{g}}{2^{g}} 
\frac{e^{2T}}{T^{2+g}}, \hbox{
\text{as} }
T \to +\infty.$$
\endproclaim

In the next section, we shall describe the necessary background
on Anosov flows, periodic orbits and homology. In section 2, we shall
prove Theorems 1 and 2.

\heading
1. Anosov Flows and Homology
\endheading

A $C^1$ flow $\phi_t : M \to M$ on a smooth compact Riemannian
manifold $M$ is called an Anosov flow if the tangent bundle admits
a continuous splitting $TM = E^0 \oplus E^s \oplus E^u$, where 
$E^0$ is the 1-dimensional bundle tangent to the flow trajectories
and where there exist constants $C>0$ and $\lambda>0$ such that
\roster
\item
$\|D\phi_t(v)\| \leq C e^{-\lambda t} \|v\|$, for all $v \in E^s$
and $t \geq 0$; and
\item
$\|D\phi_{-t}(v)\| \leq C e^{-\lambda t} \|v\|$, for all $v \in E^u$
and $t \geq 0$.
\endroster

We need some background on periodic orbits and homology 
for Anosov flows. (For more details, see the second author's 
survey in \cite{12}.)
Let $\Cal M_\phi$ denote the set of all $\phi_t$-invariant probability
measures on $M$ and, for $\mu \in \Cal M_\phi$,
let $\Phi_\mu \in H_1(M,\Bbb R)$ denote the associated
winding cycle, defined by the duality
$$\langle \Phi_\mu,[\omega] \rangle = \int \omega(\Cal X) d\mu,$$
where $[\omega]$ is the de Rham cohomology class of a closed 1-form
$\omega$ and $\Cal X$ is the vector field generating $\phi_t$. Write 
$\Cal B_\phi = \{\Phi_\mu \hbox{ : } \mu \in \Cal
M_\phi\}\subset H_1(M,\Bbb R).$
 (For geodesic flows, $\Cal B_\phi$ is
the unit ball for the Gromov-Federer stable norm on homology \cite{13}.)
The identification $H_1(M,\Bbb R)\cong\Bbb{R}^k$ defines a topology on
$H_1(M,\Bbb R)$ by considering the standard topology on $\Bbb{R}^k$
and this induces a topology on $\Cal B_\phi$ also.

Let $\mu_0 \in \Cal M_\phi$ denote the measure of
maximal entropy for $\phi_t$, i.e., the unique $\mu_0 \in
\Cal M_\phi$ for which the measure theoretic entropy $h_\phi(\mu_0)$ is
equal to the topological entropy $h$, and
write $\Phi_0 = \Phi_{\mu_0}$; this winding cycle will play a 
particularly important role. It is a fact that $\Phi_0$ lies in the
interior of $\Cal B_\phi$.

There is a (real analytic) entropy function
$\frak h : \text{\rom{int}}(\Cal B_\phi) \to \Bbb R$
defined by
$$\frak h(\rho) = \sup\{h_\phi(\mu) \hbox{ : } \Phi_\mu = \rho\}.$$
In view of the variational principle
$h = \sup\{h_\phi(\mu) \hbox{ : } \mu \in \Cal M_\phi\}$, 
$\frak h(\Phi_0) = h$ and if $\rho \neq
\Phi_0$ then $\frak h(\rho) < h$;
in particular, $\nabla \frak h(\Phi_0) =0$. In fact, it is a well-known
result that $\frak h$ is strictly concave
and that
$\Cal H = -\nabla^2\frak h(\Phi_0)$ is positive definite. 
Define a norm $\|\cdot\|$ on $H_1(M,\Bbb R) \cong \Bbb R^k$ by
$\|\rho\|^2 = \langle \rho, \Cal H \rho \rangle$. 
In particular, 
$$
\frak h(\Phi_0 + \rho) = h - \|\rho\|^2/2 +O(\|\rho\|^3).
\eqno (3)$$
when $\|\rho\|$ is sufficiently small. Also
define $\sigma >0$  by $\sigma^{-2k} = \det \Cal H$. We note that
since $H_1(M,\Bbb{R})$ has finite dimension as a real vectorspace
the norm
$\|\cdot\|$ induces the same topology as the one previously
considered. 

Let $\frak p : H^1(M,\Bbb R) \cong \Bbb R^k$ 
be the Legendre conjugate of $-\frak h$. (This
may be defined directly as a pressure function
by the formula 
$\frak p([\omega]) = P(\omega(\Cal X)) 
=\sup\{h_\phi(\mu) + \langle \Phi_\mu,[\omega] \rangle \hbox{ : } \mu \in
\Cal M_\phi\}$.)
In particular, set $\xi(\rho) = (\nabla \frak p)^{-1}(\rho)$.
Then $\xi(\Phi_0)=0$. 

\demo{Remark}
The above analysis only applies
directly when $\phi$ is a $C^{1+\epsilon}$ flow, so that the
functions $\omega(\Cal X)$ are H\"older continuous.
For the modifications we required for a flow which is only $C^1$,
see \cite{3}.
\enddemo

As in the introduction, for $\alpha \in H_1(M,\Bbb Z)/\text{torsion}$, we write
$\pi(T,\alpha) = \#\{\gamma \hbox{ : } l(\gamma) \leq T, \
[\gamma]=\alpha\}$. Now, however, we shall allow $\alpha$ to depend on
$T$ (in a linear way). To continue to take values in
$H_1(M,\Bbb Z)$, we shall define an ``integer
part'' on $H_1(M,\Bbb R)$. Choose a fundamental domain $\Cal F$ for
$H_1(M,\Bbb Z)/\text{torsion}$ as a lattice inside $H_1(M,\Bbb R)$. Then,
for $\rho \in H_1(M,\Bbb R)$, define
$\lfloor \rho \rfloor \in H_1(M,\Bbb Z)$ by $\rho - \lfloor \rho
\rfloor \in \Cal F$.

\proclaim{Proposition 1 \cite{2},\cite{9},\cite{10}}
Let $\phi_t : M \to M$ be a weak-mixing transitive Anosov flow.
If $\rho \in \text{\rom{int}}(\Cal B_\phi)$ and $\alpha_0 \in H_1(M,\Bbb
Z)/\text{\rom{torsion}}$
then 
$$\pi(T,\alpha_0+ \lfloor \rho T \rfloor) 
\sim 
C(\rho)
e^{-\langle \xi(\rho),\alpha_0 \rangle} 
e^{\langle \xi(\rho), T\rho-\lfloor T\rho \rfloor \rangle}
\frac{e^{\frak h(\rho)T}}{T^{k/2+1}}, \quad \text{\rom{ as }} T \to
+\infty,$$
where $C(\rho) = 
(\det \nabla^2 \frak h(\rho))^{1/2}/((2\pi)^{k/2}\frak h(\rho))>0$,
uniformly for $\rho$ in compact subsets of $\text{\rom{int}}(\Cal B_\phi)$.
\endproclaim

To put this in context, let us consider a fixed homology class $\alpha$.
Suppose first that
$0 \in \text{\rom{int}}(\Cal B_\phi)$; then
$$\pi(T,\alpha) \sim C(0) \frac{e^{\frak h(0)T}}{T^{1+k/2}}, \quad
\text{as } T \to +\infty$$
\cite{20}. On the other hand, if $0 \notin \Cal B_\phi$ then $\phi_t$ has
a global cross section and there are at most finitely may orbits in each
fixed class \cite{4}. If $0 \in \partial \Cal B_\phi$, the situation
is not well understood and the growth of $\pi(T,\alpha)$ may be
polynomial \cite{1} or exponential.
Regardless of these considerations, Proposition 1 gives a universal asymptotic
formula for the number of periodic orbits in homology classes
which grow like $\Phi_0 T$. 
To simplify notation, we write
$$\widetilde \pi_\alpha(T) = \pi(T,\alpha+\lfloor \Phi_0 T \rfloor).$$
We have the following corollaries of Proposition 1: 
\proclaim{Corollary 1}
For $\delta>0$ sufficiently small,
$$\lim_{T \to +\infty} \sup_{\|\alpha\| \leq \delta
T}
\left|\frac{T^{k/2+1}\widetilde \pi_\alpha(T)
e^{\langle \xi(\Phi_0+\alpha/T), T\Phi_0-\lfloor T\Phi_0 \rfloor
\rangle}}
{C(\Phi_0+\alpha/T) e^{\frak
h(\Phi_0+\alpha/T)T}}
-1\right|=0.$$
\endproclaim
This follows from Proposition 1 by using uniformity when setting
$\alpha_0=0$ and $\rho=\Phi_0+\alpha/T$. Since $\Phi_0$ is an 
interior point of $\Cal
B_\phi$ such $\rho$'s are in a compact subset of $\Cal B_\phi$ for $\delta$ 
sufficiently small.  
The following version of the
Central Limit Theorem also holds. 

\proclaim{Corollary 2}
For a Jordan set $B \subset \Bbb R^k$ whose boundary has zero measure,
$$\lim_{T \to + \infty}
\frac{1}{\pi(T)} 
\#\left\{\gamma \hbox{ : } l(\gamma) \leq T, \ \frac{[\gamma]-\lfloor\Phi_0
T\rfloor}
{\sqrt T} \in
B
\right\}
=
\frac{1}{(2\pi)^{k/2}\sigma^k}
\int_B e^{-\|x\|^2/2} dx$$
\endproclaim

This is straightforward to derived from Lemma 1 below, which in turn 
follows from
Corollary 1.

\heading
2. Proof of Theorems 1 and 2
\endheading

We now proceed to the proof of Theorem 1.
Our argument will be based on the simple yet powerful observation that
equation (2) may be replaced by
$$\pi_2^\beta(T) = \sum_{\alpha \in \Bbb Z^k} \widetilde \pi_\alpha(T)
\widetilde \pi_{\alpha + \beta}(T) \eqno (4)$$
and the properties of $\widetilde \pi_\alpha(T)$ contained in
Corollaries 1 and 2.
In particular, we shall use Corollary 1 to understand $\widetilde
\pi_\alpha(T)$ for $\|\alpha\| = O(\sqrt T)$ and Corollary 2
to show that the remaining terms make a negligible contribution.

Our first lemma, shows that in the range $\|\alpha\| = O(\sqrt T)$,
$\widetilde{\pi}_\alpha(T)$ is well approximated by a simpler function 
than the one given in Corollary 1.

\proclaim{Lemma 1} For any $\Delta >0$,
$$\sup_{\|\alpha\| \leq \Delta \sqrt T}
\left|\frac{h T\widetilde \pi_\alpha(T)}{ e^{hT}}
-\frac{e^{-\|\alpha\|^2/2T}}
{(2\pi)^{k/2}\sigma^kT^{k/2}}\right|=o\left(\frac{1}{T^{k/2}}\right).
$$
\endproclaim

\demo{Proof}
Provided $T$ is sufficiently large, $\Delta \sqrt T \leq \delta T$,
so it follows from Corollary 1 that
$$\sup_{\|\alpha\| \leq \Delta \sqrt T}
\left|\frac{T\widetilde \pi_\alpha(T)
e^{\langle \xi(\Phi_0+\alpha/T), T\Phi_0-\lfloor T\Phi_0 \rfloor
\rangle}}
{C(\Phi_0 +\alpha/T) e^{\frak
h(\Phi_0+\alpha/T)T}}
-\frac{1}{T^{k/2}}\right|=o\left(\frac{1}{T^{k/2}}\right).$$

We have $e^{\langle \xi(\Phi_0+\alpha/T), T\Phi_0-\lfloor T\Phi_0
\rfloor \rangle}=1+O(T^{-1/2})$ when $\|\alpha\| \leq \Delta \sqrt T$
so 

$$\sup_{\|\alpha\| \leq \Delta \sqrt T}
\left|\frac{T\widetilde \pi_\alpha(T)}{C(\Phi_0 +\alpha/T) e^{\frak
h(\Phi_0+\alpha/T)T}}
-\frac{1}{T^{k/2}}\right|=o\left(\frac{1}{T^{k/2}}\right).$$

Note that $C(\Phi_0)= ((2\pi)^{k/2}\sigma^k h)^{-1}$.
Since the entropy function $\frak h$ is real analytic, we have 
for $\|\alpha\| \leq \Delta \sqrt T$, 
\roster
\item"{(i)}"
$|C(\Phi_0+\alpha/T)-C(\Phi_0)|=
O(T^{-1/2})$, and, using (3),
\item"{(ii)}"
$\frak h(\Phi_0+\alpha/T)T = hT - \|\alpha\|^2/2T + O(T^{-1/2})$,
\endroster
we may replace this by
$$\sup_{\|\alpha\| \leq \Delta \sqrt T}
\left|\frac{h T\widetilde \pi_\alpha(T)}{ e^{hT}}
-\frac{e^{-\|\alpha\|^2/2T}e^{q(\alpha,T)}}
{(2\pi)^{k/2}\sigma^kT^{k/2}}\right|=o\left(\frac{1}{T^{k/2}}\right),$$
where $e^{q(\alpha,T)} \in (e^{-cT^{-1/2}},e^{cT^{-1/2}})$,
for some $c>0$. The result follows by using that $e^{q(\alpha,T)}
=O(T^{-1/2})$. \qed
\enddemo

We may then use Lemma 1 to find good approximations for 
$\sum\widetilde \pi_\alpha(T) 
\widetilde \pi_{\alpha +
\beta}(T)$ where the sum is over ${\|\alpha\| \leq \Delta \sqrt T}$

\proclaim{Lemma 2} For any $\Delta>0$,
$$\lim_{T \to +\infty} \sum_{\|\alpha\| \leq \Delta \sqrt T}
\left(\frac{(2\pi)^k \sigma^{2k} h^2T^{2+k/2}\widetilde \pi_\alpha(T) 
\widetilde \pi_{\alpha +
\beta}(T)}
{e^{2hT}} - \frac{e^{-\|\alpha\|^2/2T}e^{-\|\alpha+\beta\|^2/2T}}
{T^{k/2}}
\right) =0.$$
\endproclaim

\demo{Proof}
To shorten some of our formulae, we shall write 
$e_T(\alpha)=e^{-\|\alpha\|^2/2T}$.
We have 
$$\eqalign{
&\left|\frac{T^{2+k/2}\widetilde \pi_\alpha(T) \widetilde \pi_{\alpha +
\beta}(T)}
{C(\Phi_0)^2 e^{2hT}} - \frac{e_T(\alpha)e_T(\alpha+\beta)}
{T^{k/2}}\right| \cr
\leq &\left|\frac{T^{2+k/2}\widetilde \pi_\alpha(T) \widetilde \pi_{\alpha +
\beta}(T)}
{C(\Phi_0)^2 e^{2hT}} - \frac{Te_T(\alpha) \widetilde \pi_{\alpha
+\beta}(T)}
{C(\Phi_0)e^{hT}}\right| 
+\left|\frac{Te_T(\alpha) \widetilde \pi_{\alpha
+\beta}(T)}
{C(\Phi_0) e^{hT}} - \frac{e_T(\alpha)e_T(\alpha+\beta)}
{T^{k/2}}\right|. \cr
}$$
Applying Lemma 1, the terms on the Right Hand Side satisfy the estimates
$$o\left(\frac{T\widetilde \pi_{\alpha + \beta}(T)}{e^{hT}}\right) 
= o\left(\frac{1}{T^{k/2}}\right) \quad
\text{and} \quad
o\left(\frac{e_T(\alpha)}{T^{k/2}}\right) = 
o\left(\frac{1}{T^{k/2}}\right),$$
respectively, uniformly for $\|\alpha\| \leq \Delta \sqrt T$.
Summing over $\|\alpha\| \leq \Delta \sqrt T$
gives the result. \qed
\enddemo

Note that,
given $\epsilon >0$, it is possible to choose $\Delta>0$ sufficiently
large that
$$\frac{1}{(2\pi)^{k/2}\sigma^k}
\int_{\|x\| >\Delta} e^{-\|x\|^2/2} dx <\epsilon. \eqno (5)
$$

From Lemma 2 it is clear that we need to understand the asymptotic
behaviour of 
$$\sum_{\|\alpha\| \leq\Delta \sqrt T}
e^{-\|\alpha\|^2/2T}e^{-\|\alpha+\beta\|^2/2T}.$$ This behaviour is 
found in the next lemma:

\proclaim{Lemma 3} Given $\epsilon >0$, provided $\Delta$ is
sufficiently 
large we have
$$\pi^{k/2}\sigma^{k}(1-\epsilon) \leq \lim_{T \to +\infty} 
\frac{1}{T^{k/2}}
\sum_{\|\alpha\| \leq\Delta \sqrt T}
e^{-\|\alpha\|^2/2T}e^{-\|\alpha+\beta\|^2/2T}
 \leq \pi^{k/2}\sigma^{k}(1+\epsilon).$$
\endproclaim

\demo{Proof}
Note that
$$\eqalign{\sum_{\|\alpha\| \leq\Delta \sqrt T}
e^{-\|\alpha\|^2/2T}e^{-\|\alpha+\beta\|^2/2T}
 &= \sum_{\|\alpha\| \leq\Delta \sqrt T}
e^{-\|\alpha\|^2/T}
e^{ -(2\langle \alpha ,\Cal H \beta \rangle/T + \|\beta\|^2)/2T} \cr
&= \sum_{\|\alpha\| \leq\Delta \sqrt T}
e^{-\|\alpha\|^2/T} + O\left(\frac{1}{\sqrt T}\right).\cr
}$$
Since
$\int_{\Bbb R^k} e^{-\langle x,\Cal H x \rangle} dx 
= \pi^{k/2}/\sqrt{\det \Cal H},$
applying Lemma 2 of \cite{3} or the proof of Lemma 2.10 in \cite{16} gives
$$\lim_{T \to +\infty} \frac{1}{\pi^{k/2} \sigma^k T^{k/2}}
\sum_{\alpha \in \Bbb Z^k} e^{-\|\alpha\|^2/T} = 1.$$
Choosing $\Delta$ sufficiently large that (5) is satisfied
(and since $e^{-\|x\|^2} \leq e^{-\|x\|^2/2}$) we also have
$$\lim_{T \to +\infty}\frac{1}{\pi^{k/2} \sigma^k T^{k/2}}
\sum_{\|\alpha\| > \Delta \sqrt T} e^{-\|\alpha\|^2/T} 
=\frac{1}{\pi^{k/2} \sigma^k} 
\int_{\|x\| > \Delta} e^{-\|x\|^2} dx <\epsilon. \eqno \qed$$
\enddemo

In order to complete the proof we need a uniform upper bound on
$\widetilde \pi_\alpha(T)$ in the range where Proposition 1 gives no
information. This is provided by the following lemma:

\proclaim{Lemma 4} There exists $B>0$ such that
$$\widetilde \pi_\alpha(T) \leq B\frac{e^{hT}}{T^{1+k/2}}$$
for all $\alpha \in \Bbb Z^k$ and $T>0$.
\endproclaim

\demo{Proof}
By Corollary 1, if we fix $\delta>0$ sufficiently small then there 
exists $T_0 >0$
such that, for $T \geq T_0$ and $\|\alpha\| \leq \delta T$,
$$\widetilde \pi_\alpha(T) \leq \frac{2C(\Phi_0+\alpha/T)}
{e^{\langle \xi(\Phi_0 +\alpha/T), T\Phi_0 - \lfloor T\Phi_0 \rfloor \rangle}}
\frac{e^{\frak h(\Phi_0+\alpha/T)T}}{T^{1+k/2}} 
\leq B_0 \frac{e^{hT}}{T^{1+k/2}},$$
where
$$B_0 = 2 \sup 
 \left\{\frac{C(\Phi_0+\rho)}
{e^{\langle \xi(\Phi_0 +\rho), \rho' \rangle}} 
\hbox{ : } \|\rho\| \leq \delta, \
 \rho' \in \Cal F\right\}.$$
 
 To obtain the bound for $\|\alpha\| > \delta T$ we use large deviations
 theory. For a periodic orbit $\gamma$, let $\mu_\gamma$ denote
 the normalized Lebesgue measure around $\gamma$, i.e.,
 $$\int f \ d\mu_\gamma = \int_0^{l(\gamma)} f(\phi_t x_\gamma) dt,$$
 for any $x_\gamma \in \gamma$.  We may choose closed 1-forms
 $\omega_1,\ldots,\omega_k$ such that
 $$\frac{[\gamma]}{l(\gamma)} = 
\left(\int \omega_1(\Cal X) \ d\mu_\gamma,\ldots,
 \int \omega_k(\Cal X) \ d\mu_\gamma
 \right). \eqno(6)$$
 
 Define a set $\Cal K \subset \Cal M_\phi$ by
 $$\Cal K = \left\{\mu \in \Cal M_\phi \hbox{ : }
 \left\|\left(\int \omega_1(\Cal X) \ d\mu_\gamma,\ldots,
 \int \omega_k(\Cal X) \ d\mu_\gamma
 \right)- \Phi_0\right\| \geq \frac{\delta}{2}\right\};$$
 this is weak$^*$ compact. 
By Theorem 2.1 of \cite{8},
 $$\limsup_{T \to +\infty} 
\frac{1}{T} \log \#\{\gamma \hbox{ : } l(\gamma) \leq T, \
 \mu_\gamma \in \Cal K\} \leq h_{\Cal K} :=\sup_{\mu \in \Cal K} h_\phi(\mu).$$
 Furthermore, since $\mu_0 \notin \Cal K$, $h_{\Cal K} < h$.

Recall that $\Cal F$ is a fundamental domain for $H_1(M,\Bbb
Z)/\text{torsion}$ in $H_1(M,\Bbb R)$ and let $D$ denote
its diameter with respect to $\|\cdot\|$.
Choose $0 < \theta < 1$ 
and note that
$$\sum_{\|\alpha\| > \delta T} \widetilde \pi_\alpha(T)
= \#\{ \gamma \hbox{ : } \theta T < l(\gamma) \leq T,
\|[\gamma]-\lfloor T\Phi_0 \rfloor\| > \delta T\}
+O(e^{\theta hT}), \eqno(7)$$
where the implied constant depends only on $\theta$.

Now consider $\gamma$ with
$\theta T < l(\gamma) \leq T$. 
\comment
It is an easy observation that there exists $M>0$ such that
$\|[\gamma]\| \leq Ml(\gamma)$, for every periodic orbits $\gamma$.
\endcomment
Then
$\|[\gamma] - \lfloor T\Phi_0 \rfloor\| > \delta T$ implies that
$$\eqalign{\left\|\frac{[\gamma]}{l(\gamma)} - \Phi_0\right\|
&\geq \left\|\frac{[\gamma] - \lfloor T \Phi_0
\rfloor}{l(\gamma)}\right\|
- \left\|\frac{\lfloor T \Phi_0 \rfloor}{l(\gamma)} 
- \frac{T \Phi_0}{l(\gamma)}\right\| 
- \left\|\frac{T \Phi_0}{l(\gamma)} - \Phi_0\right\| \cr
&> \delta - \frac{D}{\theta T}- (\theta^{-1}-1)\|\Phi_0\|.  
\cr } \eqno(8)$$
If we choose $\theta$ sufficiently close to $1$ and $T_1 >0$
sufficiently large then we may assume that, provided $T \geq T_1$,
$$ \delta - \frac{D}{\theta T} -(\theta^{-1}-1)\|\Phi_0\|\geq \frac{\delta}{2}. \eqno(9)$$
Combining (8) and (9), we obtain the estimate
$$\#\{ \gamma \hbox{ : } \theta T < l(\gamma) \leq T,
\|[\gamma]-\lfloor T \Phi_0 \rfloor\| > \delta T\}
\leq 
 \#\{\gamma \hbox{ : } l(\gamma) \leq T, \
 \mu_\gamma \in \Cal K\}.$$
Applying this to (7), there exists $T_2$ such that, for $T \geq T_2$ and 
$\|\alpha\|>\delta T$,
 $$\widetilde \pi_\alpha(T) \leq e^{ h_{\Cal K} + \epsilon}.$$
 Increasing $T_2$ if necessary, we may also suppose
 that $e^{ h_{\Cal K} + \epsilon} \leq B_0 e^{hT}/T^{1+k/2}$.

Finally, we may choose $B_1>0$ so large such that, for 
$T \leq \max\{T_0,T_1,T_2\}$
 and any $\alpha \in \Bbb Z^k$,
 $\widetilde \pi_\alpha(T) \leq B_1 e^{hT}/T^{1+k/2}$.
 The proposition is thus proved with $B = \max\{B_0,B_1\}$.
\qed

\enddemo

\comment
\proclaim{Lemma 4} Given any $\delta >0$ there exist a $h_\delta<h$
and $B>0$ such that 
such that 
$$\widetilde \pi_\alpha(T)\leq B e^{h_\delta T}$$
uniformly for $\|\alpha\|>\delta T$.
\endproclaim
\demo{Proof}
PLEASE WRITE DOWN A PROOF OR AN APPROPRIATE REFERENCE TO KIFER OR
POLLICOTT. I still think that the more natural (at least in my
opinion)

$$ \pi_\alpha(T)\leq B e^{hT}/T^{k/2+1}$$
should hold uniformly for \it all \rm $\alpha$, but I cannot prove it in this
generality. In constant negative curvature it is easy but in the
weak-mixing case I do
not know. This is the (minimal?) statement that is really needed later
also. \qed
\enddemo
\endcomment

We now combine the preceding lemmas with Corollary 2 to prove Theorem 1.

\demo{Proof of Theorem 1}
Given $\epsilon>0$, choose $\Delta>0$ so that (5) is satisfied.
Consider the sum in equation (4). Lemmas 2 and 3 tell us what happens when this
sum is restricted to $\|\alpha\| \leq \Delta \sqrt T$:
we need to consider the remaining terms.
\comment
Clearly,
$$\sum_{\|\alpha\| > \Delta \sqrt T} 
\frac{T^{2+k/2}\widetilde \pi_\alpha(T) \widetilde \pi_{\alpha + \beta}(T)}
{C(\Phi_0)^2e^{2hT}} \leq 
\left(\sum_{\|\alpha\| > \Delta \sqrt T} \frac{T\widetilde \pi_\alpha(T)}
{C(\Phi_0)e^{hT}}\right)
\sup_{\|\alpha\| > \Delta \sqrt T}\left\{ \frac{T^{1+k/2} \widetilde 
\pi_{\alpha+\beta}(T)}
{C(\Phi_0)e^{hT}}\right\}.$$
Fix $\delta >0$. From Lemma 4 we conclude that when $\|\alpha\| >\delta T$, 
$\widetilde \pi_{\alpha+\beta}(T)$
has an exponential growth rate smaller than $h$. Using the uniformity
in Proposition 1, we can control also $\Delta\sqrt{T}<\|\alpha\|<\delta T$ provided that $\delta$ is sufficiently
small, so we may conclude that
$$A=\limsup_{T \to +\infty}
\sup_{\|\alpha\| > \Delta \sqrt T} 
\left\{\frac{T^{1+k/2} \widetilde \pi_{\alpha+\beta}(T)}
{C(\Phi_0)e^{hT}}\right\} < +\infty.$$
\endcomment
By Lemma 4,
$$\sum_{\|\alpha\| > \Delta \sqrt T} 
\frac{T^{2+k/2}\widetilde \pi_\alpha(T) \widetilde \pi_{\alpha + \beta}(T)}
{C(\Phi_0)^2e^{2hT}} \leq 
\frac{B}{C(\Phi_0)}
\sum_{\|\alpha\| > \Delta \sqrt T} \frac{T\widetilde \pi_\alpha(T)}
{C(\Phi_0)e^{hT}}.$$
Thus, by Corollary 2,
$$\limsup_{T \to +\infty} 
\sum_{\|\alpha\|>\Delta\sqrt{T}}\frac{T^{2+k/2}\widetilde \pi_\alpha(T) \widetilde \pi_{\alpha + \beta}(T)}
{C(\Phi_0)^2e^{2hT}} 
\leq \frac{B}{C(\Phi_0)} \left(\int_{\|x\| >\Delta} e^{-\|x\|^2/2} dx\right)  <
\frac{B}{C(\Phi_0)} \epsilon .$$
By the above estimate and Lemmas 2 and 3,
$$\eqalign{
\pi^{k/2} \sigma^k (1-\epsilon) &<
\liminf_{T \to +\infty} \sum_{\alpha \in \Bbb Z^k} 
\frac{T^{2+k/2}\widetilde \pi_\alpha(T) \widetilde \pi_{\alpha +
\beta}(T)}
{C(\Phi_0)^2e^{2hT}} \cr
&\leq
\limsup_{T \to +\infty} \sum_{\alpha \in \Bbb Z^k} 
\frac{T^{2+k/2}\widetilde \pi_\alpha(T) \widetilde \pi_{\alpha +
\beta}(T)}
{C(\Phi_0)^2e^{2hT}} < \pi^{k/2} \sigma^k (1+\epsilon) + \frac{B}{C(\Phi_0)}\epsilon. \cr }
$$
Since $\epsilon>0$ is arbitrary, this completes the proof with
$$\Cal C(\phi) = C(\Phi_0)^2 \pi^{k/2} \sigma^k 
= \frac{1}{2^k \pi^{k/2} \sigma^k h^2}. \eqno \qed$$
\enddemo

\demo{Remark}
It would be interesting to have a version of Theorem 1 where
the asymptotic behaviour was uniform in $\beta$.
A
slightly more careful version of our analysis shows that 
uniformity holds in the range $\|\beta\|  =o(\sqrt T)$ but this
is insufficient for most applications.
To obtain a stronger result,
one would need a deeper analysis of the sum
$$\sum_{\alpha \in \Bbb Z^k} e^{-(\|\alpha\|^2 +  \|\alpha+\beta\|^2)/2T}.$$
\enddemo

We conclude by proving Theorem 2.

\demo{Proof of Theorem 2}
All we need to do is to check that the constant $(g -1)^g/2^g$ is correct.
For a compact surface of constant curvature $-1$ and genus $g$,  $h=1$ and
$$\frac{1}{(2\pi)^g\sigma^{2g}}=C(\Phi_0) = (g-1)^g$$
\cite{17}, so that in this case,
$$\Cal C(\phi) = C(\Phi_0)^2 \pi^g \sigma^{2g}=
(g-1)^{2g} \pi^g \sigma^{2g}
=\frac{(g-1)^{g} }{2^{g}},$$
as required. \qed
\enddemo

\enddocument